\documentclass[11pt,a4paper,twoside,draft]{article}
\usepackage{amsthm,amsfonts,amsmath,amscd,amssymb}
\usepackage{latexsym}
\usepackage{euscript}
\usepackage{enumitem}
\usepackage{graphicx}
\usepackage{tikz}
\usepackage{caption}
\usepackage{subcaption}
\usepackage{cite}
\usepackage[utf8]{inputenc}
\usepackage[english]{babel}
\usepackage{xcolor}
\usepackage{textcomp}

\begin{document}

\begin{center}
{\bf I.P. Il'inskaya}\footnote{Kharkiv National V.N.~Karazin University, Kharkiv, Ukraine; \\ e-mail: iljinskayaip@gmail.com }
\end{center}

\begin{center}
{\bf\large Description of the class of probability measures on the group ${\mathbb Z}_2^3$ that have a trivial equivalence class}
\end{center}
\bigskip

A complete description is given of the class of probability measures on the group ${\mathbb Z}_2^3$, which are uniquely determined by the modulus of their characteristic function up to a shift.
\bigskip

\bigskip

\bigskip

{\bf 1. Introduction}
\bigskip

The interest in the question of which probability measures on $\mathbb{R}^n$ or, more generally, on a locally compact abelian group can be restored by modulus of its characteristic function up to a shift and central symmetry, arose in physics (see \cite{bib1}, \cite{bib2}, and the references therein). To formulate the problem, let us introduce the required notation and definitions. Let $X$ be a locally compact abelian group, $Y$ be its character group, $(x,y)$ be the value of character $y\in Y$ on the element $x\in X$, $M^1(X)$ be the set of probability measures on the group $X$.
Characteristic function of the measure $\mu\in M^1(X)$ is determined by the formula
\begin{equation}\label{eq1}
\widehat{\mu}(y)=\int_X(x,y)\mu(dx),\quad y\in Y\,.
\end{equation}
\medskip

{\bf Definition 1}.
We say that the measures $\mu,\nu\in M^1(X)$ are {\it equivalent} and write $\mu\sim\nu$ if
\begin{equation}\label{eq2}
|\widehat{\nu}(y)|=|\widehat{\mu}(y)|,\quad y\in Y\,.
\end{equation}
\medskip

Let $\mu_x$ be the shift of measure $\mu$ by element $x\in X$ and let $\mu^-$ be the measure obtained from $\mu$ by central symmetry:
\begin{equation*}
\mu_x(E)=\mu(E+x),\quad \mu^-(E)=\mu(-E)\,,
\end{equation*}
where $E$ is a Borel set in $X$. It is easy to see that $\mu\sim\mu_x$, $\mu\sim\mu^-$, $\mu\sim\mu_x^-$ for any measure $\mu\in M^1(X)$ and any $x\in X$.
\medskip

{\bf Definition 2}. We say that a measure $\mu\in M^1(X)$ has {\it a trivial equivalence class} if only the measures $\mu_x$ and $\mu_x^-$, $x\in X$, are equivalent to it.
\medskip

The set of measures from $M^1(X)$ that have a trivial equivalence class is denoted by $TEC(X)$. Note that the normal distribution on $\mathbb{R}$ belongs to the class $TEC(\mathbb{R})$ (this follows from Cramer's theorem on the decompositions of the normal distribution) and the Poisson distribution does not belong to the class $TEC(\mathbb{R})$.

Papers \cite{bib3} -- \cite{bib10} are devoted to the task of a complete or partial description of the class $TEC(X)$ for some groups $X$. In papers
\cite{bib4}, \cite{bib5}, the question of the triviality of the equivalence class for uniform distributions on intervals of the group $\mathbb{Z}_n$ of residue classes modulo $n$, on the Cartesian product of intervals of the group $\mathbb{Z}^l$, on the unit ball in $\mathbb{R}^l$ is studied.
In \cite{bib10} a criterion was obtained for a two-point measure on $\mathbb{Z}_n$ to have a trivial equivalence class.
In \cite{bib4}, a necessary and sufficient condition is found for the fact that a generalized Poisson distribution on the group
$\mathbb{Z}_2^l$, $l\geq 2$, whose spectral measure is proportional to the Haar measure, belongs to the class
$TEC(\mathbb{Z}_2^l)$ (Theorem~1 below). In \cite{bib8}, a necessary and sufficient condition was obtained for a generalized Poisson distribution with an arbitrary spectral measure to belong to the class $TEC(\mathbb{Z}_2^2)$. In \cite{bib9}, a criterion is obtained for the fact that a generalized Poisson distribution on the group $\mathbb{Z}_2^3$, whose spectral measure is arbitrarily distributed on any three generators of the group, belongs to the class $TEC(\mathbb{Z}_2^3)$ (Theorem~2 below).

Sometimes it is possible to obtain a complete description of the class $TEC(X)$. It is easy to see that $TEC(\mathbb{Z}_2)=M^1(\mathbb{Z}_2)$. Classes $TEC(\mathbb{Z}_3)$ and $TEC(\mathbb{Z}_4)$ are fully described in \cite{bib10}, class $TEC(\mathbb{Z}_2^2)$ is described in \cite{bib8} (Theorem~3 below). In this paper, we obtain a complete description of the class $TEC(\mathbb{Z}_2^3)$ (Theorem~4 below).

For comparison we formulate some of the results from the papers \cite{bib4}, \cite{bib8}, \cite{bib9}.
We need the following definition.
\medskip

{\bf Definition 3}.
Let $\rho$ be a finite measure on the group $X$. {\it The generalized Poisson distribution} with spectral measure $\rho$ is the distribution
$$
\Pi_{\rho}=e^{-\rho(X)}\sum_{k=0}^{\infty}\frac{1}{k!}\rho^{\ast\, k}\,.
$$
\medskip

The characteristic function of the distribution $\Pi_{\rho}$ has the form
$$
\widehat{\Pi}_{\rho}(y)=\exp\Big\{\int_X[(x,y)-1]\rho(dx)\Big\}\,.
$$
Let $m_X$ denote the Haar measure on the group $X$. For a compact group $X$ we put $m_X(X)=1$. The following theorem gives a necessary and sufficient condition for the measure $\Pi_{\lambda m}$ (here $m$ is the Haar measure on the group $\mathbb{Z}_2^l$, $\lambda>0$)
to belong to the class $TEC(\mathbb{Z}_2^l)$, $l=2, 3, \ldots$.
\medskip

{\bf Theorem 1} ([2]). {\it The measure $\Pi_{\lambda m}$ belongs to the class $TEC(\mathbb{Z}_2^l)$, $l=2, 3, \ldots$, if and only if $\lambda<\ln 3$.}
\medskip

The following Theorem~2 gives a necessary and sufficient condition that the generalized Poisson distribution on
the group $\mathbb{Z}_2^3$,   whose spectral measure is arbitrarily distributed on any three generators of the group, belongs to the class  $TEC(\mathbb{Z}_2^3)$.
\medskip

{\bf Theorem 2} ([7]). {\it
Let $\pi$ be a generalized Poisson distribution on the group $\mathbb{Z}_2^3$, whose spectral measure is concentrated on any three generators and assigns masses $a$, $b$, $c$ to them. A measure $\pi$  belongs to the class $TEC(\mathbb{Z}_2^3)$ if and only if the system of inequalities
\begin{equation*}
\left\{
\begin{aligned}
&e^{-2a}+e^{-2b}+e^{-2(a+b)}>1\,,\\
&e^{-2b}+e^{-2c}+e^{-2(b+c)}>1\,,\\
&e^{-2a}-e^{-2b}+e^{-2c}+e^{-2(a+b)}+e^{-2(b+c)}-e^{-2(c+a)}+e^{-2(a+b+c)}>1\\
\end{aligned}
\right.
\end{equation*}
is satisfied, or one of the two systems of inequalities obtained from this system by cyclic permutations of the parameters $a$, $b$, $c$ is satisfied.}
\medskip

The following theorem gives the full description of the class $TEC(\mathbb{Z}_2^2)$. Let $\mu\in M^1(\mathbb{Z}_2^2)$.
We denote $a_{\text{max}}=\max\{\mu(\{x\}):x\in \mathbb{Z}_2^2\}$. Let $S(\mu)$ be the support of the measure $\mu$, $|C|$ be the number of elements of the set $C$.
\medskip

{\bf Theorem 3} ([6]). {\it
Class $TEC(\mathbb{Z}_2^2)$ contains the following measures and only them:

$1)$ all measures $\mu$ for which $|S(\mu)|\leq 2$;

$2)$ all measures $\mu$ for which $|S(\mu)|=3$ and $a_{\max}\geq 1/2$;

$3)$ all measures $\mu$ for which $|S(\mu)|=4$ and one of the following  two conditions is satisfied:

\hspace{0.31cm} a$)$ $a_{\max}>1/2$,

\hspace{0.31cm} b$)$ $a_{\max}<1/2$ and the sum of masses of some two elements of the group is equal to the sum of the masses of the other two elements.}
\medskip

We note that $\mu\notin TEC(\mathbb{Z}_2^2)$ if $|S(\mu)|=4$ and $a_{\text{max}}=1/2$.
\medskip
\medskip

{\bf 2. Main result}
\bigskip

To formulate the main result, we need some notation.

Let $\mathbb{Z}_2=\{0, 1\}$ be additive group of the residue classes modulo $2$,
$$
X=\mathbb{Z}_2^3=\left\{x=(\alpha, \beta, \gamma): \alpha, \beta, \gamma\in \{0, 1\}\right\}\,.
$$
We denote the elements of the group $\mathbb{Z}_2^3$ as follows:
$$
\begin{aligned}
&x_0=(0,0,0),\quad x_1=(1,0,0),\quad x_2=(0,1,0),\quad x_3=(0,0,1)\,,\\
&x_4=(1,1,0),\quad x_5=(1,0,1),\quad x_6=(0,1,1),\quad x_7=(1,1,1)\,.\\
\end{aligned}
$$
For the zero element $x_0$ of the group $\mathbb{Z}_2^3$ we will also use the notation $0$.

The group of characters of the group $\mathbb{Z}_2^3$ is isomorphic to $\mathbb{Z}_2^3$. The value of the character
$y=(\xi,\eta,\zeta)$, $\xi, \eta, \zeta\in \{0, 1\}$ on the element $x=(\alpha, \beta ,\gamma)\in \mathbb{Z}_2^3$ is determined by the formula  $$
(x,y)=(-1)^{\alpha\xi+\beta\eta+\gamma\zeta}\,.
$$

Let us denote by $\mathfrak{A}_1$ the set of all subgroups of the group $\mathbb{Z}_2^3$ isomorphic to $\mathbb{Z}_2$:
$$
\mathfrak{A}_1=\left\{H_i=\{x_0,x_i\}:i=1,\ldots ,7\right\}\,.
$$
Let us denote by $\mathfrak{A}_2$ the set of all subgroups of the group $\mathbb{Z}_2^3$ isomorphic to $\mathbb{Z}_2^2$:
$$
\mathfrak{A}_2=\left\{K_i:i=1,\ldots ,7\right\}\,,
$$
\begin{equation}\label{eq3}
\begin{aligned}
&K_1=\{x_0,x_2,x_3,x_6\},\quad K_2=\{x_0,x_1,x_3,x_5\},\quad K_3=\{x_0,x_1,x_2,x_4\},\\
&K_4=\{x_0,x_1,x_6,x_7\},\quad K_5=\{x_0,x_3,x_4,x_7\},\quad K_6=\{x_0,x_2,x_5,x_7\},\\
&K_7=\{x_0,x_4,x_5,x_6\}\,.\\
\end{aligned}
\end{equation}
In what follows, the representation of the group $\mathbb{Z}_2^3$ as a direct sum of a group from
$\mathfrak{A}_1$ and a group from $\mathfrak{A}_2$ will play an important role.

Let $\mu\in M^1(\mathbb{Z}_2^3)$. The values of the masses of the measure $\mu$ on the elements of the group $\mathbb{Z}_2^3$ will be denoted by
$$
a_i=\mu(\{x_i\}),\quad i=0,1,\ldots,7\,.
$$
Thus we have $a_i\geq 0$, $\sum_{i=0}^7a_i=1$. We denote
$$
a_{\max}:=\max\{a_i: i=0,1,\ldots,7\}\,.
$$

From (\ref{eq1}) we find the general form of the characteristic function $\widehat{\mu}$, $\mu\in M^1(\mathbb{Z}_2^3)$:
\begin{equation}\label{eq4}
\widehat{\mu}(y)=\sum\limits_{i=0}^7a_i(x_i,y)\,,\qquad a_i\geq 0, \,\,\, \sum_{i=0}^7a_i=1\,.
\end{equation}

Since all non-zero elements of the group $\mathbb{Z}_2^3$ have order~2, the central symmetry is the identity mapping. Therefore Definition~2 that measure $\mu\in M^1(\mathbb{Z}_2^3)$ has a trivial equivalence class is simplified:
$$
\mu\in TEC(\mathbb{Z}_2^3) \,\,\, \text{if and only if its shifts are equivalent to it}.
$$
Let $E$ be a subset of $\mathbb{Z}_2^3$. We denote
$$
u(E):=\max\left\{\mu(\{x\}): x\in E\right\}\,,\quad v(E):=\min\left\{\mu(\{x\}): x\in E\right\}\,.
$$
\medskip

{\bf Definition 4}.
Let us denote by $U(E)$ the set of measures $\mu\in M^1(\mathbb{Z}_2^3)$ satisfying the condition
$$
2u(E)>\mu(E)\,.
$$
\medskip
In other words, $U(E)$ is the set of such measures for which the maximum mass of the elements of the set $E$ is greater than the sum of the remaining masses of the elements of this set.

\medskip

{\bf Definition 5}.
Let us denote by
$V(E)$ the set of measures $\mu\in M^1(\mathbb{Z}_2^3)$ satisfying the condition
$$
1/2+2v(E)<\mu(E)\,.
$$
\medskip

In other words, $V(E)$ is the set of such measures for which the sum of $1/2$ and minimum mass of the elements of the set $E$ is smaller than the sum of the remaining masses of the elements of this set.

Let $K\in \mathfrak{A}_2$. The coset of the group $\mathbb{Z}_2^3$ in the subgroup $K$, different from $K$, will be denoted by $\overline{K}$.

For convenience of references, we formulate the following statement in the form of a lemma. It easy follows from the description of subgroups of
$\mathbb{Z}_2^3$.
\medskip

{\bf Lemma 1}. {\it
Let $H\in \mathfrak{A}_1$, $K\in \mathfrak{A}_2$, and $H\cap K=\{0\}$. Then$:$

$1)$ there are exactly three subgroups $L_i\in \mathfrak{A}_2$ containing $H$$;$

$2)$ non-zero elements of the subgroup $K$ belong to the different subgroups $L_i$, $i=1,2,3$$;$

$3)$ there are exactly four subgroups $K^{(j)}\in \mathfrak{A}_2$, $j=1, 2, 3, 4$, such that $H\cap K^{(j)}=\{0\}$$;$

$4)$ for each $j_0=1, 2, 3, 4$ and for each $j\not=j_0$, there is a unique element $z_j\in \overline{K^{(j)}}$
such that $z_j\not\in H$, $z_j\not\in K^{(j_0)}$.}
\medskip

We illustrate Lemma~1 with the help of Figures~1 and~2. If we take
$H=\{0,x_1\}$, then $L_1=\{0, x_1,x_2,x_4\}$, $L_2=\{0, x_1, x_3, x_5\}$, $L_3=\{0, x_1, x_6, x_7\}$
(see Figure~2, where subgroups $L_1$,  $L_2$, $L_3 $ are marked with different shadings).
Let $K^{(j_0)}=\{0,x_2,x_3,x_6\}$ and $K^{(j)}=\{0,x_3,x_4,x_7\}$ (see Figure~1). Then $\overline{K^{(j)}}=\{x_1,x_2,x_5,x_6\}$, and since the conditions $z_j\notin H$, $z_j\notin{K^{(1)}}$ must be satisfied, we see that $z_j\not=x_1$, $z_j\not=x_2$, $z_j\not=x_6$. Therefore $z_j=x_5$.


\begin{picture}(360,160)
\put(70,70){\circle*{3}}
\put(55,70){$x_0$}
\put(35,35){\circle*{3}}
\put(20,35){$x_2$}
\put(140,70){\circle*{3}}
\put(143,75){$x_3$}
\put(70,140){\circle*{3}}
\put(56,143){$x_1$}
\put(105,35){\circle*{3}}
\put(110,33){$x_6$}
\put(35,105){\circle*{3}}
\put(20,100){$x_4$}
\put(140,140){\circle*{3}}
\put(143,143){$x_5$}
\put(105,105){\circle*{3}}
\put(110,103){$x_7$}
\put(35,35){\line(1,0){70}}
\put(70,70){\line(1,0){70}}
\put(35,105){\line(1,0){70}}
\put(70,140){\line(1,0){70}}
\put(35,35){\line(0,1){70}}
\put(105,35){\line(0,1){70}}
\put(70,70){\line(0,1){70}}
\put(140,70){\line(0,1){70}}
\put(35,35){\line(1,1){35}}
\put(35,105){\line(1,1){35}}
\put(105,35){\line(1,1){35}}
\put(105,105){\line(1,1){35}}
\put(70,70){\line(-1,1){35}}
\put(140,70){\line(-1,1){35}}

\put(35,35){\vector(-1,-1){14}}
\put(26,19){$\beta$}
\put(140,70){\vector(1,0){14}}
\put(156,65){$\gamma$}
\put(70,140){\vector(0,1){14}}
\put(73,149){$\alpha$}
\multiput(45,40)(5,0){13}{\circle*{1}}
\multiput(50,45)(5,0){13}{\circle*{1}}
\multiput(55,50)(5,0){13}{\circle*{1}}
\multiput(60,55)(5,0){13}{\circle*{1}}
\multiput(65,60)(5,0){13}{\circle*{1}}
\multiput(70,65)(5,0){13}{\circle*{1}}
\multiput(132,73)(-7,0){9}{\circle*{1}}
\multiput(123,82)(-7,0){9}{\circle*{1}}
\multiput(114,91)(-7,0){9}{\circle*{1}}
\multiput(105,100)(-7,0){9}{\circle*{1}}
\put(80,5){Fig. 1}
\put(270,70){\circle*{3}}
\put(257,75){$x_0$}
\put(235,35){\circle*{3}}
\put(220,35){$x_2$}
\put(340,70){\circle*{3}}
\put(343,75){$x_3$}
\put(270,140){\circle*{3}}
\put(256,143){$x_1$}
\put(305,35){\circle*{3}}
\put(310,33){$x_6$}
\put(235,105){\circle*{3}}
\put(220,100){$x_4$}
\put(340,140){\circle*{3}}
\put(343,143){$x_5$}
\put(305,105){\circle*{3}}
\put(310,100){$x_7$}
\put(235,35){\line(1,0){70}}
\put(270,70){\line(1,0){70}}
\put(235,105){\line(1,0){70}}
\put(270,140){\line(1,0){70}}
\put(235,35){\line(0,1){70}}
\put(305,35){\line(0,1){70}}
\put(270,70){\line(0,1){70}}
\put(340,70){\line(0,1){70}}
\put(235,35){\line(1,1){35}}
\put(235,105){\line(1,1){35}}
\put(305,35){\line(1,1){35}}
\put(305,105){\line(1,1){35}}
\put(305,35){\line(-1,1){35}}
\put(305,105){\line(-1,1){35}}
\put(235,35){\vector(-1,-1){14}}
\put(226,19){$\beta$}
\put(340,70){\vector(1,0){14}}
\put(356,65){$\gamma$}
\put(270,140){\vector(0,1){14}}
\put(273,149){$\alpha$}
\multiput(275,75)(0,5){13}{\circle*{1}}
\multiput(280,75)(0,5){13}{\circle*{1}}
\multiput(285,75)(0,5){13}{\circle*{1}}
\multiput(290,75)(0,5){13}{\circle*{1}}
\multiput(295,75)(0,5){13}{\circle*{1}}
\multiput(300,75)(0,5){13}{\circle*{1}}
\multiput(305,75)(0,5){13}{\circle*{1}}
\multiput(310,75)(0,5){13}{\circle*{1}}
\multiput(315,75)(0,5){13}{\circle*{1}}
\multiput(320,75)(0,5){13}{\circle*{1}}
\multiput(325,75)(0,5){13}{\circle*{1}}
\multiput(330,75)(0,5){13}{\circle*{1}}
\multiput(335,75)(0,5){13}{\circle*{1}}
\multiput(272,73)(5,-5){7}{\line(0,1){62}}
\multiput(240,45)(7,7){5}{\circle*{1}}
\multiput(240,52)(7,7){4}{\circle*{1}}
\multiput(240,59)(7,7){3}{\circle*{1}}
\multiput(240,66)(7,7){5}{\circle*{1}}
\multiput(240,73)(7,7){5}{\circle*{1}}
\multiput(240,80)(7,7){5}{\circle*{1}}
\multiput(240,87)(7,7){5}{\circle*{1}}
\multiput(240,94)(7,7){5}{\circle*{1}}
\multiput(240,101)(7,7){5}{\circle*{1}}
\put(280,5){Fig. 2}
\end{picture}

To introduce Definition~6, we need the following notation.
Let $H=\{0,g\}\in\mathfrak{A}_1$, $K\in\mathfrak{A}_2$, and $K\cap H=\{0\}$. Let $L_1,L_2,L_3$ be all subgroups from $\mathfrak{A}_2$ containing $H$ (see items~1 and~2 of Lemma~1).
For all $i=1,2,3$ we denote by $t_i$ an element from $L_i$ that belongs to $\overline K$ and $t_i\neq g$.
Let $s_i\in L_i$ and $s_i\neq 0,g,t_i$. Thus, $L_i=\{0,g,s_i,t_i\}$, $i=1,2,3$, $K=\{0,s_1,s_2,s_3\}$.
\medskip

{\bf Definition 6}.
We denote by $W(H,K)$ the set of measures $\mu\in M^1(\mathbb{Z}_2^3)$ such that the following three equalities hold
$$
\mu(\{0, t_i\})=\mu(\{g, s_i\})\,,\quad i=1,2,3\,.
$$
\medskip

For example, if $H=\{x_0,x_1\}$, $K=\{x_0,x_2,x_3,x_6\}$
(see Figure~2), then $L_1=\{x_0,x_1,x_2,x_4\}$, $L_2=\{x_0,x_1,x_3,x_5\}$, $L_3=\{x_0,x_1,x_6,x_7\}$,
$t_1=x_4$, $t_2=x_5$, $t_3=x_7$, $s_1=x_2$, $s_2=x_3$, $s_3=x_6$.
In this case three equalities from the Definition~6 have the form
$$
\mu(\{x_0,x_4\})=\mu(\{x_1,x_2\})\,,\mu(\{x_0,x_5\})=\mu(\{x_1,x_3\})\,,\mu(\{x_0,x_7\})=\mu(\{x_1,x_6\})\,.
$$

The main result of the paper is the following theorem. Note that item I.2(c) of this theorem uses the items~3 and~4 of Lemma~1. Item~II.2 uses item~3 of this lemma. Without loss of generality we may assume that the condition $a_{\max}=a_0$ is satisfied.
\medskip

{\bf Theorem 4.} {\it Let $\mu\in M^1(\mathbb{Z}_2^3)$ and $a_{\max}=a_0$. Then $\mu\in TEC(\mathbb{Z}_2^3)$ if and only if there exists a decomposition of the group $\mathbb{Z}_2^3$ $\,:$
\begin{equation}\label{eq5}
\mathbb{Z}_2^3=X_1\oplus X_2\,,\quad X_1\in\mathfrak{A}_1\,\,\, X_2\in\mathfrak{A}_2\,,
\end{equation}
such that conditions ${\rm I}.1$ and ${\rm I}.2$ are satisfied for $a_{\max}\leq 1/4$ and one of the conditions $\rm{II}.1 - \rm{II}.4$ is satisfied for $a_{\max}>1/4$\,$:$

${\rm\bf I}.{\bf 1}$. The projection of the measure $\mu$ on $X_2$ parallel to $X_1$ is equal to the Haar measure $m_{X_2}$ of the subgroup $X_2$.

${\rm\bf I}.{\bf 2}$. At least one of the following three requirements is true$\,:$

\hspace{0.2cm} a$)$ The projection of the measure $\mu$ on $X_1$ parallel to $X_2$  is equal to the Haar measure $m_{X_1}$ of the subgroup $X_1\,;$

\hspace{0.2cm} b$)$ The sum of the masses of some two elements of the subgroup $X_2$ is equal to the sum of the masses of its other two elements$\,;$

\hspace{0.2cm} c$)$ Let $K^{(j)}$, $j=1,2,3,4$, be all subgroups of $\mathfrak{A}_2$ that do not contain $X_1$. There is $j_0$ such that the condition $\mu\in U(K^{(j_0)})$ is satisfied and the equality
$$
u\Big(\overline{K^{(j)}}\Big)=\mu(\{z_j\})
$$
holds for any $j\not=j_0$, where $z_j\in\overline{K^{(j)}}$,  $z_j\notin X_1$,  $z_j\notin K^{(j_0)}$.

${\rm\bf II}.{\bf 1}$. For any subgroup $K\in\mathfrak{A}_2$ the following conditions are true$\,:$

\hspace{0.2cm} a$)$ $\mu\in U(K)$ or $\mu\in U(\overline{K})$$\,;$

\hspace{0.2cm} b$)$ $\mu\in V(K)$ or $\mu\in V(\overline{K})$.

${\rm\bf II}.{\bf 2}$. The following two conditions hold$\,:$

\hspace{0.2cm}  a$)$ $\mu\in W(X_1,K)$ for at least one subgroup $K\in\mathfrak{A}_2$ that does not contain $X_1$$\,;$

\hspace{0.2cm}  b$)$ $\mu\in U(K)$ or $\mu\in U(\overline{K})$ for any subgroup $K\in\mathfrak{A}_2$ that does not contain $X_1$.

${\rm\bf II}.{\bf 3}$. The following two conditions hold$\,:$

\hspace{0.2cm} a$)$ There is a two-point subset $E\subset X_2$ such that for both elements $g\in X_1$ the equality is valid
\begin{equation}\label{eq6}
\mu(g+E)=\mu(g+(X_2\setminus E))\,;
\end{equation}

\hspace{0.2cm} b$)$ $\mu\in U(K)$ or $\mu\in U(\overline{K})$ for each subgroup  $K\in\mathfrak{A}_2$, $K\not=X_2$.

${\rm\bf II}.{\bf 4}$. The following three conditions hold$\,:$

\hspace{0.2cm} a$)$ The projection of the measure $\mu$ on $X_1$ parallel to $X_2$ is equal to the Haar measure $m_{X_1}$ of the group $X_1$;

\hspace{0.2cm} b$)$ $\mu\in V(K)$ or $\mu\in V(\overline{K})$ for any subgroup $K\in\mathfrak{A}_2$, $K\not=X_2$$\,;$

\hspace{0.2cm} c$)$ $\mu\in U(K)$ or $\mu\in U(\overline{K})$ for any subgroup $K\in\mathfrak{A}_2$.
}

\medskip

Note that in the statement of condition II.1 of Theorem~4, decomposition~(\ref{eq5}) is not used.

Let us indicate a simple sufficient condition for a measure to belong to the class $TEC(\mathbb{Z}_2^3)$.
\medskip

{\bf Corollary 1.} {\it If  $a_{\max}>5/6$, then $\mu\in TEC(\mathbb{Z}_2^3)$.}
\medskip

{\bf Proof.} Let $K$ be an arbitrary subgroup from $\mathfrak{A}_2$. Since, without loss of generality, we can assume that $a_{\max}=a_0$, we have that condition $\mu\in U(K)$, and hence condition II.1($a$), are satisfied.
Since $v(K)<1/6$ and $\mu(K)>5/6$, we see that
$$
1/2+2v(K)<1/2+2\cdot 1/6=5/6<\mu(K)\,.
$$
Therefore, condition II.1($b$) is satisfied.
Hence, $\mu\in TEC(\mathbb{Z}_2^3)$.
\bigskip

{\bf 3. Examples}
\bigskip

In this section we give the examples of measures that satisfy various conditions of Theorem~4. In examples~1~---~7, we will assume that the subgroups $X_1$ and $X_2$ appearing in the formulation of Theorem~4 are as follows:
$$
X_1=\{x_0, x_1\}, \quad X_2=\{x_0, x_2, x_3, x_6\}.
$$

Let us give an example of a distribution that satisfies conditions~I.1 and~I.2($a$).
\medskip

{\bf Example 1}.
We consider the distribution with masses
\begin{align*}
&a_0=1/8+4\varepsilon, & &a_1=1/8-4\varepsilon, & &a_2=1/8-3\varepsilon, & &a_3=1/8-2\varepsilon,\\
&a_4=1/8+3\varepsilon, & &a_5=1/8+2\varepsilon, & &a_6=1/8+\varepsilon, & &a_7=1/8-\varepsilon\\
\end{align*}
for $0\leq\varepsilon\leq 1/32$. (Note that for $\varepsilon=0$ we obtain the Haar measure on the group $\mathbb{Z}_2^3$.)
Since $a_0+a_1=a_2+a_4=a_3+a_5=a_6+a_7=1/4$, we see that condition~I.1 is satisfied. Since
$a_0+a_2+a_3+a_6=1/2$, we see that condition~I.2($a$) is satisfied.

Let us give an example of a distribution that does not satisfy condition~I.2($a$), but satisfies conditions~I.1 and I.2($b$).
\medskip

{\bf Example 2}. We consider the distribution with masses
\begin{equation*}
\begin{aligned}
&a_0=a_3=1/4-\varepsilon,\,\,\, a_1=a_5=\varepsilon,\,\,\, a_2=a_6=2\varepsilon\,,\\
&a_4=a_7=1/4-2\varepsilon \quad (0<\varepsilon<1/8).
\end{aligned}
\end{equation*}
Since $a_0+a_1=a_2+a_4=a_3+a_5=a_6+a_7=1/4$ and $a_0+a_2=a_3+a_6$, we see that I.1 and I.2($b$) are true. Since
$a_0+a_2+a_3+a_6=1/2+2\varepsilon$, condition~I.2($a$) is not satisfied.

Let us give an example of a distribution that satisfies conditions~I.1 and~I.2($c$).

\medskip

{\bf Example 3}. We consider the distribution $\mu$ with masses
\begin{equation*}
\begin{aligned}
&a_0=1/4-\varepsilon,\,\,\, a_1=\varepsilon,\,\,\, a_2=a_3=a_6=2\varepsilon\,,\\
&a_4=a_5=a_7=1/4-2\varepsilon \quad (0<\varepsilon<1/28).
\end{aligned}
\end{equation*}
Since $a_0+a_1=a_2+a_4=a_3+a_5=a_6+a_7=1/4$, we obtain that condition~I.1 is satisfied. Let us show that condition~I.2($c$) is satisfied. All subgroups of
$\mathfrak{A}_2$ that do not contain  subgroup $X_1$ (see Lemma~1, item~3) are as follows:
\begin{equation*}
\begin{aligned}
&K^{(1)}=K_1=\{x_0,x_2,x_3,x_6\},\quad K^{(2)}=K_5=\{x_0,x_3,x_4,x_7\},\\
&K^{(3)}=K_6=\{x_0,x_2,x_5,x_7\},\quad K^{(4)}=K_7=\{x_0,x_4,x_5,x_6\}\,.\\
\end{aligned}
\end{equation*}
Let us take $j_0=1$. It is easy to see that
$$
\max\{a_i:x_i\in K^{(1)}\}=a_0=1/4-\varepsilon>6\varepsilon=a_2+a_3+a_6
$$
for $0<\varepsilon<1/28$.
Therefore, $\mu\in U(K^{(1)})$.

For $j=2$ we have
$$
\max\{a_i:x_i\in \overline{K^{(2)}}\}=\max\{a_1,a_2,a_5,a_6\}=a_5=1/4-2\varepsilon\,\,\,(0<\varepsilon\leq 1/16),
$$
and we can take $z_2=x_5$ (see Lemma~1, item~4):
$z_2=x_5\notin X_1$, $z_2\notin K^{(1)}$.

For $j=3$ we have
$$
\max\{a_i:x_i\in\overline{K^{(3)}}\}=\max\{a_1,a_3,a_4,a_6\}=a_4=1/4-2\varepsilon \,\,\, (0<\varepsilon\leq 1/16),
$$
and we can take $z_3=x_4$:
$z_3=x_4\notin X_1$, $z_3\notin K^{(1)}$.

For $j=4$ we have
$$
\max\{a_i:x_i\in\overline{K^{(4)}}\}=\max\{a_1,a_2,a_3,a_7\}=a_7=1/4-2\varepsilon \,\,\, (0<\varepsilon\leq 1/16),
$$
and we can take $z_4=x_7$:
$z_4=x_7\notin X_1$, $z_4\notin K^{(1)}$.

Thus, condition~I.2($c$) is satisfied. Therefore, $\mu\in TEC(\mathbb{Z}_2^3)$ for $0<\varepsilon<1/28$.

\medskip

{\bf Example 4}.
Condition II.1 is satisfied for any measure for which  $a_{\text{max}}>5/6$.
This follows from Corollary~1 of Theorem~4.

Let us give an example of a measure that satisfies condition~II.2.

\medskip

{\bf Example 5}. We consider the distribution $\mu$ with masses
\begin{equation*}
\begin{aligned}
&a_0=8/24+3\varepsilon,\,\,\, a_1=7/24+3\varepsilon\,,\\
&a_2=a_3=a_6=2/24-\varepsilon,\,\,\, a_4=a_5=a_7=1/24-\varepsilon\,,
\end{aligned}
\end{equation*}
where $0\leq\varepsilon\leq 1/24$.
Let $K=X_2=\{x_0, x_2, x_3, x_6\}$. The subgroups containing the subgroup $X_1$ (see Lemma~1, item~1) are as follows:
$$
L_1=\{x_0, x_1,x_2,x_4\},\,\,\, L_2=\{x_0, x_1, x_3, x_5\},\,\,\, L_3=\{x_0, x_1, x_6, x_7\}.
$$
Notice that $x_2\in L_1, x_3\in L_2, x_6\in L_3$. We put $s_1=x_2$, $s_2=x_3$, $s_3=x_6$.
Then $t_1=x_4$, $t_2=x_5$, $t_3=x_7$. It is easy to see that condition $\mu\in W(X_1,K)$ is satisfied, since
\begin{equation*}
\begin{aligned}
&\mu(\{x_0,x_4\})=\mu(\{x_1,x_2\})=9/24+2\varepsilon,\\
&\mu(\{x_0,x_5\})=\mu(\{x_1,x_3\})=9/24+2\varepsilon,\\
&\mu(\{x_0,x_7\})=\mu(\{x_1,x_6\})=9/24+2\varepsilon.
\end{aligned}
\end{equation*}
Therefore, condition~II.2($a$) is satisfied. Let us check the fulfillment of condition~II.2($b$).
For the subgroup  $\{x_0, x_2, x_3, x_6\}\not\supset X_1$ we have
$$
\max\{a_0,a_2,a_3,a_6\}=a_0=8/24+3\varepsilon>3\cdot (2/24-\varepsilon)=a_2+a_3+a_6.
$$
For the subgroup $\{x_0, x_3, x_4, x_7\}\not\supset X_1$ we have
\begin{equation*}
\begin{aligned}
\max\{a_0,a_3,a_4,a_7\}&=a_0=8/24+3\varepsilon>4/24-3\varepsilon=(2/24-\varepsilon)+2\cdot (1/24-\varepsilon)\\
                       &=a_3+a_4+a_7.
\end{aligned}
\end{equation*}
Similar inequalities are also valid for subgroups $\{x_0, x_2, x_5, x_7\}\not\supset X_1$,
$\{x_0, x_4, x_5, x_6\}\not\supset X_1$. Therefore, condition~II.2($b$) is satisfied, and hence $\mu\in TEC(\mathbb{Z}_2^3)$.

Let us give an example of a distribution that satisfies condition~II.3.

\medskip

{\bf Example 6}. We consider a distribution $\mu$ with masses
\begin{equation*}
\begin{aligned}
&a_0=1/4+\varepsilon,\,\,\, a_1=a_4=a_5=a_7=1/8-\varepsilon,\\
&a_2=a_3=1/8+\varepsilon,\,\,\, a_6=\varepsilon \qquad (1/16<\varepsilon<1/8)\,.
\end{aligned}
\end{equation*}
We take
$E=\{x_2,x_3\}$. Then $X_2\setminus E=\{x_0,x_6\}$.
We check the fulfillment of condition~II.3($a$).

Let us show that equality (\ref{eq6}) holds for $g=x_0$:
\begin{equation*}
\begin{aligned}
&\mu(\{x_2,x_3\})=2\cdot(1/8+\varepsilon)=1/4+2\varepsilon,\\
&\mu(\{x_0,x_6\})=1/4+\varepsilon+\varepsilon=1/4+2\varepsilon.
\end{aligned}
\end{equation*}

Let us show that equality (\ref{eq6}) holds for $g=x_1$:
\begin{equation*}
\begin{aligned}
&\mu(\{x_1+\{x_2,x_3\}\})=\mu(\{x_4,x_5\})=2\cdot(1/8-\varepsilon)=1/4-2\varepsilon,\\
&\mu(\{x_1+\{x_0,x_6\}\})=\mu(\{x_1,x_7\})=2\cdot(1/8-\varepsilon)=1/4-2\varepsilon.\\
\end{aligned}
\end{equation*}

Let us check that condition~II.3($b$) is satisfied. For $K=\{x_0, x_1, x_3, x_5\}$, condition $\mu\in U(K)$ means that
$$
a_1+a_3+a_5=(1/8-\varepsilon)+(1/8+\varepsilon)+(1/8-\varepsilon)=3/8-\varepsilon<1/4+\varepsilon=a_0,
$$
and this inequality is satisfied for any $\varepsilon$ such that $1/16<\varepsilon<1/8$. The condition
$\mu\in U(K)$ for $K=\{x_0, x_1, x_2, x_4\}$, $\{x_0, x_3, x_4, x_7\}$, $\{x_0, x_2, x_5, x_7\}$ can be checked in the same way.
For $K=\{x_0, x_4, x_5, x_6\}$, condition $\mu\in U(K)$ means that
$$
a_4+a_5+a_6=2\cdot(1/8-\varepsilon)+\varepsilon=1/4-\varepsilon<1/4+\varepsilon=a_0.
$$
This inequality is true for any $\varepsilon>0$. The condition $\mu\in U(K)$ for $K=\{x_0, x_1, x_6, x_7\}$ can be verified in a similar way.
Thus, condition~II.3($b$) is verified.
Therefore, $\mu\in TEC(\mathbb{Z}_2^3)$.

Let us give an example of a distribution that satisfies condition~II.4.

\medskip

{\bf Example 7}. We consider the distribution $\mu$ with masses
$$
a_0=1/2-\varepsilon,\,\,\, a_2=a_3=a_6=\varepsilon/3,\,\,\, a_1=a_4=a_5=a_7=1/8\,\,\,\,\,  (0<\varepsilon<3/16).
$$
Since $a_0+a_2+a_3+a_6=1/2$, condition~II.4($a$) is fulfilled. Let us check that condition~II.4($b$) is satisfied.
Let $K=K_2=\{x_0, x_1, x_3, x_5\}$. We have:
\begin{equation*}
\begin{aligned}
&a_0+a_1+a_3+a_5=(1/2-\varepsilon)+1/8+\varepsilon/3+1/8=3/4-2\varepsilon/3,\\
&1/2+2\min\{a_0,a_1,a_3,a_5\}=1/2+2\varepsilon/3.
\end{aligned}
\end{equation*}
Condition $\mu\in V(K)$ takes the form $1/2+2\varepsilon/3<3/4-2\varepsilon/3$ and is satisfied for $0<\varepsilon<3/16$.
For other subgroups $K_i$, $i=3,\ldots,7$, the check is the same.

Let us check the fulfillment of condition~II.4($c$). If $K=K_1=\{x_0, x_2, x_3, x_6\}$, we have
$$
a_0=1/2-\varepsilon>3\cdot\varepsilon/3=\varepsilon=a_2+a_3+a_6 \,\,\, (0<\varepsilon<1/4).
$$
If $K=K_2=\{x_0, x_1, x_3, x_5\}$, we have
$$
a_1+a_3+a_5=1/8+\varepsilon/3+1/8=1/4+\varepsilon/3<1/2-\varepsilon=a_0\,\,\, (0<\varepsilon<3/16).
$$
For other subgroups $K_i$, $i=3,\ldots,7$, the check is the same as for $K_2$.
Therefore, $\mu\in TEC(\mathbb{Z}_2^3)$ for $0<\varepsilon<3/16$.

Let us give two examples of the measures that do not belong to the class $TEC(\mathbb{Z}_2^3)$. Let $X_1$ and $X_2$ be the arbitrary subgroups of the group $\mathbb{Z}_2^3$ that satisfy the condition~(\ref{eq5}).
\medskip

{\bf Example 8}.
We consider the distribution $\mu$ with masses
$$
a_0=1/4+\varepsilon,\,\,\, a_1=a_2=\ldots =a_7=(3/4-\varepsilon)/7  \quad (0<\varepsilon\leq 1/20).
$$
For this distribution we have $a_{\max}> 1/4$.
Let us show that conditions~II.1($a$), II.2($b$), II.3($b$), II.4($c$) are not satisfied.
Let $K$ be an arbitrary subgroup of $\mathfrak{A}_2$, $K=\{x_0,x_i,x_j,x_k\}$,
where $i,j,k$ are different and nonzero.
Condition $\mu\notin U(K)$ means that $a_0\leq a_i+a_j+a_k$,
that is
$$
1/4+\varepsilon\leq 3\cdot(1/7)(3/4-\varepsilon)\,.
$$
This inequality holds for $\varepsilon\leq 1/20$.
Condition $\mu\notin U(\overline{K})$
means that
$$
(1/7)(3/4-\varepsilon)\leq 3\cdot(1/7)(3/4-\varepsilon)\,.
$$
This inequality holds for $\varepsilon\leq 3/4$.
Thus, $\mu\notin TEC(\mathbb{Z}_2^3)$.

In the previous example, the measure satisfied condition $a_{\max}>1/4$.
In the following example, the condition $a_{\max}\leq 1/4$ is valid.
\medskip

{\bf Example 9}. We consider the distribution with masses
$$
a_0=1/4-\varepsilon,\,\,\, a_1=a_2=\ldots =a_7=(3/4+\varepsilon)/7  \quad (\varepsilon\in [0,1/8)\cup (1/8,1/4]).
$$
We show that condition~I.1 is not satisfied.
We have
$$
a_0+a_i=5/14-(6/7)\varepsilon\not=1/4\,,\quad \mbox{if}\quad \varepsilon\not=1/8\,,\quad i=1,2,\ldots ,7\,.
$$
Therefore, condition~I.1 is not satisfied. Thus, $\mu\notin TEC(\mathbb{Z}_2^3)$.
Note that this measure is the Haar measure of the group $\mathbb{Z}_2^3$ if $\varepsilon=1/8$. Therefore, it belongs to the class  $TEC(\mathbb{Z}_2^3)$.
\bigskip

{\bf 4. Derivation of Theorems 1 -- 3 from Theorem 4}
\bigskip

Let us show that Theorem~4 implies Theorem~1 in the case of $l=3$ and Theorems~2 and~3.

We show that the sufficiency in Theorem~1 for the case $l=3$ follows from item~II.1 of Theorem~4. Indeed, if $m$ is the Haar measure of the group $\mathbb{Z}_2^3$, we have $\widehat{\Pi}_{\lambda m}(y)=1$ for $y=0$ and
$\widehat{\Pi}_{\lambda m}(y)=e^{-\lambda}$ in other cases. Therefore, as it is easy to check, the masses of the measure $\Pi_{\lambda m}$ are as follows: the mass of the zero element of the group $\mathbb{Z}_2^3$ is equal to $(1+7e^{-\lambda})/8$, and the masses of the remaining elements of the group are equal to $(1-e^{-\lambda})/8$. Therefore, condition~II.1($a$) of Theorem~4 takes the form
$(1+7e^{-\lambda})/8> 3(1-e^{-\lambda})/8$ and is satisfied for $e^{-\lambda}>1/5$,
while condition~II.1($b$) takes the form $1/2<(1+7e^{-\lambda})/8+(1-e^{-\lambda})/8$
and is satisfied for $e^{-\lambda}>1/3$.
Thus, the conditions of item~II.1 of Theorem~4 are satisfied for the measure $\Pi_{\lambda m}$ for $\lambda<\ln 3$, which proves the sufficiency in Theorem~1.

The sufficiency in Theorem~2 also follows from item~II.1 of Theorem~4. If one of the three systems of inequalities from the formulation of Theorem~2 is true, the conditions of item~II.1 of Theorem~4 are satisfied. (The proof of this fact is lengthy, so we omit it.)

Let us show how the sufficiency in Theorem~3 follows from Theorem~4. To do this, we consider several cases:

1) $a_{\max}>1/2$, $|S(\mu)|=2,3,4$;

2) $a_{\max}=1/2$, $|S(\mu)|=2$;

3) $a_{\max}=1/2$, $|S(\mu)|=3$;

4) $a_{\max}<1/2$, $|S(\mu)|=4$ and the sum of some two masses of elements from $S(\mu)$ is equal to the sum of two other masses.

Item 4) is divided into three sub-items, differing in the number of coinciding  masses:

$i$) all masses of elements are equal (the Haar measure on $\mathbb{Z}_2^2$),

$ii$) among the masses, there are two pairs of the same one, but not all of these masses are equal,

$iii$) all masses are different or there are exactly two identical masses among them.

(The case when there are exactly three identical masses among them is impossible, since in this case the condition that the sum of some two masses is equal to the sum of two other masses cannot be satisfied.)

It is easy to see that condition~II.1 of Theorem~4 is satisfied under condition~1) of Theorem~3.

Let us show that, under condition~4($i$), conditions~I.1 and~I.2($a$) of Theorem~4 are satisfied. Without loss of generality, we can assume that $S(\mu)=\{x_0,x_1,x_3,x_5\}$. Then $a_0=a_1=a_3=a_5=1/4$, $a_2=a_4=a_6=a_7=0$.
Let us put $X_1=\{x_0, x_2\}$, $X_2=\{x_0, x_3, x_4, x_7\}$.
Since $a_0+a_2=a_1+a_4=a_3+a_6=a_5+a_7=1/4$, condition~I.1 of Theorem~4 is satisfied. Since $a_0+a_3+a_4+a_7=a_1+a_2+a_5+a_6=1/2$, condition~I.2($a$) is satisfied.
Therefore $\mu\in TEC(\mathbb{Z}_2^3)$.

Let us show that, under condition~4($ii$), condition~II.2 of Theorem~4 is satisfied. Without loss of generality, we can assume that
$S(\mu)=\{x_0,x_1,x_3,x_5\}$, $a_{\text{max}}=a_0<1/2$,  $a_0+a_5=a_1+a_3$, $a_0=a_1$, $a_3=a_5$, $a_0\not=a_3$.
Let $X_1=\{x_0, x_1\}$, $X_2=\{x_0, x_2, x_3, x_6\}$
(see Fig.~2). We put $K=\{x_0, x_2, x_3, x_6\}\not\supset X_1$.
Let us show that $\mu\in W(X_1,K)$.
We denote
$$
L_1=\{x_0, x_1,x_2,x_4\},\,\,\, L_2=\{x_0, x_1, x_3, x_5\},\,\,\, L_3=\{x_0, x_1, x_6, x_7\}
$$
(see item~1 of Lemma~1).
We have $L_i\supset X_1$ for $i=1,2,3$.
It follows from the conditions on the masses $a_i$ that the equalities
\begin{equation*}
\begin{aligned}
&\mu(\{x_0,x_4\})=\mu(\{x_1,x_2\}),\quad \mu(\{x_0,x_5\})=\mu(\{x_1,x_3\})\,,\\
&\mu(\{x_0,x_7\})=\mu(\{x_1,x_6\})
\end{aligned}
\end{equation*}
are valid. Therefore, condition~II.2($a$) is satisfied.
It is easy to see that condition $\mu\in U(K)$ is satisfied: $a_0>a_2+a_3+a_6=a_3$.
The fulfillment of condition~II.2($b$) for the remaining subgroups from $\mathfrak{A}_2$ that do not contain $X_1$ is verified similarly. So, $\mu\in TEC(\mathbb{Z}_2^3)$.

Let us show that, under condition~2), condition~II.2 of Theorem~4 is satisfied. Without loss of generality, we can assume that
$S(\mu)=\{x_0,x_1\}$, $a_0=a_1=1/2$.
The fulfillment of condition~II.2 is checked in exactly the same way as for the case~4($ii$), if we put $a_3=a_5=0$ in the previous argument.

Let us show that, under condition~4($iii$), condition~II.3 of Theorem~4 is satisfied. We can assume that $S(\mu)=\{x_0,x_1,x_3,x_5\}$, $a_0>a_1\geq a_3>a_5$, $a_0+a_5=a_1+a_3$.
Let us put $X_1=\{x_0, x_2\}$, $X_2=\{x_0, x_1, x_3, x_5\}$, $E=\{x_0,x_5\}$.
Then $X_2\setminus E=\{x_1, x_3\}$.
It is clear that
\begin{equation*}
\begin{aligned}
&\mu(E)=a_0+a_5=a_1+a_3=\mu(X_2\setminus E)\,,\\
&\mu(x_2+E)=\mu(\{x_2,x_7\})=0,\quad \mu(x_2+(X_2\setminus E))=\mu(\{x_4,x_6\})=0\,.
\end{aligned}
\end{equation*}
Therefore, condition~II.3($a$) is satisfied.
Let us check the fulfillment of condition~II.3($b$).
If $K=\{x_0, x_2, x_3, x_6\}$, we have $\mu\in U(K)$, since $a_0>a_2+a_3+a_6=a_3$. Similarly, condition $\mu\in U(K)$ is satisfied for all other subgroups $K\in\mathfrak{A}_2$, $K\not=X_2$.
So, condition~II.3($b$) is satisfied. Therefore, $\mu\in TEC(\mathbb{Z}_2^3)$.

Let us show that, under condition~3), condition~II.3 of Theorem~4 is also satisfied. We can assume that
$S(\mu)=\{x_0,x_1,x_3\}$, $a_0=a_1+a_3=1/2$, $a_1>0$, $a_3>0$.
Condition~II.3 is verified in exactly the same way as in the case~4($iii$), if we put $a_5=0$ in the previous reasoning.

To verify the validity of Theorem~1 (for $l=3$) and Theorems~2, 3 in the direction of necessity, it is necessary to show that if the conditions of one of them are not satisfied, then condition~I.1 or condition~I.2 of Theorem~4 is not satisfied, and also all four conditions~II.1~--~II.4 are not satisfied. Checking this is easy, but lengthy. Therefore, we omit it.

Note that the derivation of Theorem~3 from Theorem~4 turned out to be more difficult than the proof of Theorem~3 itself.
\bigskip


{\bf 5. Scheme of the proof of Theorem~4}
\bigskip

The proof of Theorem~3 ([6]) is easy. The proof of Theorem~2 ([7])  is rather laborious.
The proof of Theorem~4, based on the results of computer calculations, is much more complicated, multi-step, and very lengthy. Let us briefly describe the scheme of the proof of Theorem~4.

Let $\mu\in M^1(\mathbb{Z}_2^3)$ be a given measure, $\nu\in M^1(\mathbb{Z}_2^3)$ be a measure equivalent to it, that is, condition~(\ref{eq2}) is satisfied. The question is the following: under which conditions on the measure $\mu$ does the equality $\nu=\mu_x$ hold for some $x$? In other words, we look for the conditions on the measure $\mu$ under which   the following equality holds
\begin{equation}\label{eq7}
\widehat{\nu}(y)=(x_j,y)\widehat{\mu}(y)
\end{equation}
for some $j=0,1,\ldots,7$.
Let the characteristic function of the measure $\mu$ be of the form~(\ref{eq4}) and the characteristic function of the measure $\nu$ be equal to
\begin{equation*}
\widehat{\nu}(y)=\sum\limits_{i=0}^7b_i(x_i,y)\,.
\end{equation*}
If in equality (\ref{eq2}), which is valid for all $y=(\xi,\eta,\zeta)$, we assume
$\xi,\eta,\zeta=0,1$, then it will be written as a system of eight equalities:
\begin{equation*}
\begin{aligned}
&b_0+b_1+b_2+b_3+b_4+b_5+b_6+b_7=1\,,\\
&|b_0-b_1+b_2+b_3-b_4-b_5+b_6-b_7|=|a_0-a_1+a_2+a_3-a_4-a_5+a_6-a_7|\,,\\
&|b_0+b_1-b_2+b_3-b_4+b_5-b_6-b_7|=|a_0+a_1-a_2+a_3-a_4+a_5-a_6-a_7|\,,\\
&|b_0+b_1+b_2-b_3+b_4-b_5-b_6-b_7|=|a_0+a_1+a_2-a_3+a_4-a_5-a_6-a_7|\,,\\
&|b_0-b_1-b_2+b_3+b_4-b_5-b_6+b_7|=|a_0-a_1-a_2+a_3+a_4-a_5-a_6+a_7|\,,\\
&|b_0-b_1+b_2-b_3-b_4+b_5-b_6+b_7|=|a_0-a_1+a_2-a_3-a_4+a_5-a_6+a_7|\,,\\
&|b_0+b_1-b_2-b_3-b_4-b_5+b_6+b_7|=|a_0+a_1-a_2-a_3-a_4-a_5+a_6+a_7|\,,\\
&|b_0-b_1-b_2-b_3+b_4+b_5+b_6-b_7|=|a_0-a_1-a_2-a_3+a_4+a_5+a_6-a_7|\,.\\
\end{aligned}
\end{equation*}
Expanding the modulus in the seven equalities of this system, we see that~(\ref{eq2}) is equivalent to a set of~$2^7=128$ systems of linear equations with eight unknowns~$b_i$ and eight given parameters~$a_i$. Computer calculations of the paper~\cite{bib9} give the solution of each system. Careful analysis of the solutions shows that~128 systems are divided into four sets according to the type of solutions, we denote them by $\cal A$, $\cal B$, $\cal C$, $\cal D$.

The set $\cal A$ consists of eight systems whose solutions $b_0,b_1,\ldots,b_7$ are such that the function $\widehat{\nu}(y)$ has the form~(\ref{eq7}).
The set~$\cal B$ consists of eight systems whose solutions are such that
$$
\widehat{\nu}(y)=(x_j,y)\sum\limits_{i=0}^7(1/4-a_i)(x_i,y)\,,\quad j=0,1,\ldots,7\,.
$$

To describe systems from sets $\cal C$ and $\cal D$, we will need the following substitutions of the indices $0,1,\ldots,7$
($e$ is the unit substitution):
\begin{equation}\label{eq8}
\begin{aligned}
&\sigma_1=e,& &\sigma_2=(12)(56),& &\sigma_3=(13)(46), & &\sigma_4=(12)(37),\\
&&&\sigma_5=(24)(67),& &\sigma_6=(35)(67),& &\sigma_7=(24)(35)\,.\\
\end{aligned}
\end{equation}
Substitution $\sigma_k$, $k=2,3,\ldots,7$, transforms the subgroup $K_1\in\mathfrak{A}_2$ into the subgroup $K_k$ (see~(\ref{eq3})).

The set $\cal C$ consists of~56 systems, which are subdivided into~7 subsets ${\cal C}_k$, $k=1,2,\ldots,7$, of~8 systems each. For solutions $b_0,b_1,\ldots,b_7$ of the systems from the subset ${\cal C}_k$, the function $\widehat{\nu}(y)$ has the form
$$
\widehat{\nu}(y)=(x_j,y)\psi_k(y)\,,\quad j=0,1,\ldots,7\,,
$$
where
$$
\psi_1(y)=\sum\limits_{i=0}^7c_i(x_i,y)\,,
$$

\begin{equation}\label{eq9}
\begin{aligned}
&c_0=(-a_0+a_2+a_3+a_6)/2, & &c_4=(a_1-a_4+a_5+a_7)/2,\\
&c_1=(-a_1+a_4+a_5+a_7)/2, & &c_5=(a_1+a_4-a_5+a_7)/2,\\
&c_2=(a_0-a_2+a_3+a_6)/2,  & &c_6=(a_0+a_2+a_3-a_6)/2,\\
&c_3=(a_0+a_2-a_3+a_6)/2,  & &c_7=(a_1+a_4+a_5-a_7)/2\,.\\
\end{aligned}
\end{equation}
To obtain the coefficients of the function $\psi_k(y)$, $k=2,3,\ldots,7$, it is necessary to apply in formulas~(\ref{eq9}) the same substitution $\sigma_k$ from~(\ref{eq8}) to the coefficients $a_i$ and $c_i$.

The set $\cal D$ consists of~56 systems, which are subdivided into~7 subsets ${\cal D}_k$, $k=1,2,\ldots,7$, of~8 systems each. For solutions of systems from the subset ${\cal D}_k$, the function $\widehat{\nu}(y)$ has the form
$$
\widehat{\nu}(y)=(x_j,y)\varphi_k(y)\,,\quad j=0,1,\ldots,7\,.
$$
where
$$
\varphi_1(y)=\sum\limits_{i=0}^7d_i(x_i,y)\,,\quad d_i=1/4-c_i\,,\,\, i=0,1,\ldots,7\,.
$$
Coefficients of the function $\varphi_k(y)$, $k=2,3,\ldots,7$, are obtained from coefficients of the function $\varphi_1(y)$ by using the substitution $\sigma_k$ applied to coefficients $a_i$ and $d_i$.

A solution $b_0,b_1,\ldots,b_7$ of any of the systems is called {\it trivial\,} if the function $\widehat{\nu}(y)$ satisfies the equality~(\ref{eq7})
(this means that the measure $\nu$ is a shift of the measure $\mu$).
We call a solution {\it non-trivial} if $b_i\geq 0$ for all $i$, and equality~(\ref{eq7}) is not satisfied for any $j$. Since we are only interested in non-negative $b_i$, we will say that the system {\it has no solution} if $b_i<0$ for some $i$.

Since all systems of the set $A$ have trivial solutions, to prove Theorem~4, it is necessary to find conditions on $a_i$ under which each system from the sets $\cal B$, $\cal C$, $\cal D$ either has a trivial solution or has no solutions.

The finding of a condition under which a particular system has no solutions is not difficult. For example, at least one of the coefficients of the function $\psi_1(y)$ is negative if and only if one of the two inequalities holds:
\begin{equation*}
\begin{aligned}
&2\max\{a_0,a_2,a_3,a_6\}>a_0+a_2+a_3+a_6\,,\\
&2\max\{a_1,a_4,a_5,a_7\}>a_1+a_4+a_5+a_7\,.
\end{aligned}
\end{equation*}
It is much more difficult to find out under what conditions a particular system or group of systems has a trivial solution.
To do this, one needs to find out when the solution of this system $b_0,b_1,\ldots,b_7$ coincides with the set of coefficients of some of the~8 functions $(x_j,y)\widehat{\mu}(y)$,  $j=0,1,2,\ldots,7$. And this is the hardest part of the proof.
\medskip

Let us explain what the conditions of Theorem~4 mean. The inequality
$a_{\text{max}}\leq 1/4$
is a necessary and sufficient condition for the fact that all systems of the set $\cal B$
have solutions. Condition~I.1 is a necessary and sufficient condition that all systems of the set $\cal B$
have trivial solutions. In order for all systems of sets $\cal C$ and $\cal D$ to have trivial solutions, it is necessary and sufficient that one of the conditions I.2($a$), I.2($b$) is satisfied.
Condition II.1($a$) is a necessary and sufficient condition that all systems of the set $\cal C$ have no solutions. Condition II.1($b$) is a necessary and sufficient condition that all systems of the set $\cal D$
have no solutions. Conditions I.1 and I.2($c$) are necessary and sufficient conditions that all systems of the set
$\cal B$, part of the systems of the set $\cal C$, and part of the systems of the set $\cal D$ have trivial solutions, and the remaining systems of the sets $\cal C$ and $\cal D$ have no solutions.
Condition
$a_{\text{max}}>1/4$
and fulfillment of one of the conditions II.2, II.3 --- this is a necessary and sufficient condition that all systems of the set $\cal B$,
part of the systems of the set $\cal C$ and part of the systems of the set $\cal D$ have no solutions, and the remaining systems of the sets $\cal C$ and $\cal D$ have trivial solutions.
Conditions
$a_{\text{max}}>1/4$
and II.4 --- this is a necessary and sufficient condition that all systems of the sets $\cal B, \cal C$
and part of the systems of the set $\cal D$ have no solutions, and the remaining systems of the set $\cal D$ have trivial solutions.
Several cases not described here (for example, when all systems of the set $\cal C$ have trivial solutions, and all systems of the set $\cal D$
have no solutions) are impossible.

I would like to thank G.M.~Fel'dman for useful discussions of the manuscript.




\begin{thebibliography}{99}

\bibitem{bib3} H. Carnal, M. Dozzi, \emph{On a decomposition problem for multivariate probability measures}, J.~Multivar. Anal., {\bf 31} (1989), 165--177.
\bibitem{bib4} H. Carnal, G.M. Fel'dman, \emph{Phase retrieval for probability measures on abelian groups}, I, J.~Theor. Probab., {\bf 8}:3 (1995), 717--725.
\bibitem{bib5} H. Carnal, G.M. Fel'dman, \emph{Phase retrieval for probability measures on abelian groups}, II, J.~Theor. Probab., {\bf 10}:4 (1997), 1065--1074.
\bibitem{bib6} H. Carnal, G.M. Fel'dman, \emph{On a property of entire characteristic functions of finite order with real zeros}, DAN, {\bf 366}:2 (1999), 162--163 (in Russian).
\bibitem{bib7} H. Carnal, G.M. Feldman, \emph{A stability property for probability measures on Abelian groups},
{\it Statistics and probability Letters}, {\bf 49} (2000), 39--44.
\bibitem{bib8} I.P. Il'inskaya, \emph{Phase retrieval for probability measures on character group of the Cantor-Walsh group}, Dopovidi NAN Ukrainy, \textnumero~8 (2003), 11--14 (in Russian).
\bibitem{bib9} I.P. Il'inskaya, D.S. Neguritsa, \emph{On probability measures on the group of Walsh functions with trivial equvivalence class},  Ukrainian Mathematical Journal, {\bf 65}:5 (2013), 717--721 (in Russian).
\bibitem{bib10} I.P. Il'inskaya, \emph{Phase retrieval for probability measures on cyclic groups}, Mathematychni Studii, {\bf 46}:1 (2016), 89 -- 95.
\bibitem{bib2} N. E. Hurt, \emph{Phase retrieval and zero crossings. Mathematical methods in image
reconstruction}\,, Kluwer Academic
Publishers, Dordrecht etc., 1989, 304 pp.
\bibitem{bib1} J. Rosenblatt, \emph{Phase retrieval}, Comm. Math. Phys., {\bf 95} (1984), 317--343.

\end{thebibliography}
\end{document}